# Critical Clearing Time Sensitivity for Inequality Constrained Systems

Chetan Mishra, *Member, IEEE*, Reetam Sen Biswas, *Student Member, IEEE*, Anamitra Pal, *Member, IEEE*, and Virgilio A. Centeno, *Senior Member, IEEE*

*Abstract*—With the growth of renewable generation (RG) and the development of associated ride through curves serving as operating limits, during disturbances, on violation of these limits, the power system is at risk of losing large amounts of generation. In order to identify preventive control measures that avoid such scenarios from manifesting, the power system must be modeled as a constrained dynamical system. For such systems, the interplay of feasibility region (man-made limits) and stability region (natural dynamical system response) results in a positively invariant region in state space known as the constrained stability region (CSR). After the occurrence of a disturbance, as it is desirable for the system trajectory to lie within the CSR, critical clearing time (CCT) must be defined with respect to the CSR instead of the stability region as is done traditionally. The sensitivity of CCT to system parameters of constrained systems then becomes beneficial for planning/revising protection settings (which impact feasible region) and/or operation (which impact dynamics). In this paper, we derive the first order CCT sensitivity of generic constrained power systems using the efficient power system trajectory sensitivity computation, pioneered by Hiskens in [1]. The results are illustrated for a single-machine infinite-bus (SMIB) system as well as a multi-machine system in order to gain meaningful insight into the dependence between ability to meet constraints, system stability, and changes occurring in power system parameters, such as, mechanical power input and inertia.

## I. NOMENCLATURE

| | |
|---|---|
| $\partial A(X)$ | Stability Region Boundary of a set $X$ |
| $A(X)$ | Stability Region of a set $X$ |
| $A_c(X)$ | Constrained Stability Region of a set $X$ |
| $f(x)$ | Vector Field |
| $g(x)$ | Equality constraint |
| $h(x)$ | Inequality constraint vector function |
| $H(x)$ | Scalar function given by $\prod_{k=1}^{n_h} h(x)$ |
| $m$ | Dimension of dependent states |
| $n$ | Dimension of state space |
| $n_h$ | Number of individual inequality constraints |
| $\aleph_X$ | Connected component of a set $X$ |
| $p$ | Parameter |
| $t_{cr}$ | Critical Clearing Time (CCT) of Base Critical Trajectory |
| $t_{cl}$ | Fault Clearing Time of Generalized Fault Trajectory |
| $T$ | Maximum Time the Base Critical Post-Fault Trajectory is Simulated |
| $W^s(X)$ | Stable Manifold of a set $X$ |
| $W^u(X)$ | Unstable Manifold of a set $X$ |
| $W^c(X)$ | Center Manifold of a set $X$ |
| $W^{c-}(X)$ | Set of trajectories lying on a subset of $W^c(X)$ converging to $X$ as $t \to \infty$ |
| $W^{c+}(X)$ | Set of trajectories lying on a subset of $W^c(X)$ converging to $X$ as $t \to -\infty$ |
| $x_{cu}$ | Controlling UEP (CUEP) of Generalized Critical Trajectory |
| $x_{cu_0}^{post}$ | CUEP of Base Critical Trajectory |
| $x_{end}^{post}$ | End point of Generalized Post-Fault Trajectory |
| $x_e$ | Original System Equilibrium Point $\{x \mid f(x) = 0\}$ |
| $x_e^s$ | Original System SEP |
| $x_e^u$ | Original System UEP |
| $x_H$ | Pseudo EP or Point on Feasibility Boundary $\{x \mid H(x) = 0\}$ |
| $x_H^0$ | Semi-saddle Pseudo EPs $\{x \mid H(x) = 0, \frac{\partial H}{\partial x} \times f(x) = 0, \frac{\partial H}{\partial x} \neq 0, \frac{\partial \left(\frac{\partial H}{\partial x} \times f(x)\right)}{\partial x} \times f(x) \neq 0\}$ |
| $x_s$ | Stable Equilibrium Point (SEP) |
| $x_{s_0}^{pre}$ | SEP of Base Pre-Fault System |
| $x_H^s$ | Stable pseudo EP |
| $x$ | State vector |
| $x_{cl}$ | State vector value at $t_{cl}$ for Generalized Fault Trajectory |
| $x_u$ | Unstable Equilibrium Point (UEP) |
| $x_H^u$ | Unstable pseudo EP |
| $x_T$ | State Value at the End of Base Critical Post-Fault Trajectory |
| $y$ | Dependent state of DAE system |
| $\varphi(x_0, t, p)$ | Trajectory of a parametric dynamical system $\dot{x} = f(x, p)$ with parameter $p$ starting from the point $x_0$ |

## II. INTRODUCTION

As opposed to the traditional approach of tripping renewable generation (RG) offline during disturbances seen at the point of common coupling (PCC), RG sources are currently made to "ride through" these disturbances. This has become necessary because systems with significant RG penetration could be at risk of collapse (particularly loss of equilibrium) if a large quantity of such generation was lost at the time of need [2], [3]. However, RGs cannot be made to ride through *every possible scenario*, especially when islanding scenarios manifest. Therefore, ride through curves were devised in the form of time dependent voltage and frequency limits at the PCC of renewable generators, violation of which resulted in their tripping. Violations of ride-through or other protection settings (such as under-voltage or under-frequency load shedding) often result in undesirable changes to the system. An example is the under-voltage load shedding and widespread tripping of distributed generation due to a delayed voltage recovery post-fault. In small or weakly connected systems, such changes can cause voltage collapse or loss of synchronism. Large systems that are strongly connected and have many controllable devices, such as the Eastern Interconnection, are

Chetan Mishra (email: chetan.mishra@dominionenergy.com) is an Engineer in Electric Transmission Operations Engineering at Dominion Energy, Richmond, VA-23220, USA.

Reetam Sen Biswas (email: rsenbisw@asu.edu) is a PhD student and Anamitra Pal (email: anamitra.pal@asu.edu) is an Assistant Professor in the School of Electrical, Computer, and Energy Engineering at Arizona State University, Tempe, AZ-85287, USA.

Virgilio A. Centeno (email: virgilio@vt.edu) is an Associate Professor in the Bradley Department of Electrical and Computer Engineering at Virginia Tech, Blacksburg, VA-24061, USA.

more robust to such follow-up events. However, analyses of recent blackouts have indicated that outage of a single element at a crucial time can destabilize a large system [4]. Therefore, operators and planners must take measures to avoid situations where a disturbance can trigger a cascade. For the measures to be accurate, the power system must be modeled as a constrained dynamical system [5] [6], where the constraints being focused on are based on preference/necessity, whose violations should be prevented at all costs.

Critical clearing time (CCT) refers to the maximum time that can be taken to clear a fault and still retain system stability. For dynamically constrained systems, one must also *not* violate the constraints. This additional requirement considerably complicates the desired starting region in the state space for post-fault trajectories. Correct decisions that enhance CCT in dynamically constrained systems will help reduce the likelihood of triggering follow-up events as well as increase dynamical stability of the system. However, to make correct decisions, knowledge of the dependence between CCT and the changes occurring in the system parameters is necessary. Generating this knowledge for constrained power systems is the primary focus of this work. For example, sensitivity of CCT to Q injection at each bus for a fault that resulted in tripping of large amounts of RGs could help plan resources such as placement of STATCOMs. Another application could be to understand an approximate dependence between system inertia and likelihood of occurrence of under frequency load shedding to establish critical inertia levels [7].

In the past, brute force approaches for CCT sensitivity computation were proposed, which primarily relied on numerical integration. Ayasun [8] reduced the multi-machine system to a single machine infinite bus (SMIB) system to evaluate sensitivities; an approach that is known to have limitations for multi-machine systems [9]. Chiodo and Lauria [10] used linear regression to understand the mapping between logarithm of CCT and loading. Nguyen [11] and Laufenberg [12] computed sensitivity of angle and speed trajectory in the post-fault phase w.r.t. fault clearing time. Since marginally stable and unstable trajectories start close to each other but later grow apart (with the former converging to a stable equilibrium point (SEP) while the latter does not); transient stability will be evident in the trajectory sensitivities. Nguyen further used the sensitivity of parametric energy function value at controlling unstable equilibrium point (CUEP) and sensitivity of fault trajectory to estimate the sensitivity of CCT. The most recent relevant work in this area was by Dobson [13] where the sensitivity of stable manifold of CUEP was used in conjunction with fault on trajectory sensitivity to estimate the sensitivity of CCT to parameter changes. The derivation in that paper was for unconstrained ordinary differential equation (ODE)-governed systems. In this paper, in continuation of Dobson's work, we derive CCT sensitivities for *inequality constrained dynamical systems* (henceforth referred to as *constrained systems*).

The rest of the paper is structured as follows. In Section III. the stability theory for constrained systems along with characterization of the quasi-stability boundary is presented. Derivation for sensitivities of various critical manifolds is provided in Section IV. Extension of the approach to constrained differential algebraic equation (DAE)-governed systems is presented in Section V. A brief description of the overall computation is provided in Section VI. Finally, to gain visual insight into the problems being addressed here, the results obtained for an SMIB system and a multi-machine system are presented in Section VII.

### III. STABILITY OF CONSTRAINED SYSTEMS

A constrained system is defined by the state equation,
$$\dot{x}_{n\times 1} = f(x,p)_{n\times 1} \qquad (1)$$
$$h(x,p)_{n_h\times 1} > 0$$

The first vector equation defines the evolution of states, $x$, while the second one defines a feasibility region with the feasibility boundary given by $\{x|(\prod_k h_k(x,p)) = H(x,p) = 0\}$. The system being analyzed is parametric with parameter $p$ but the stability properties is discussed for a fixed value of $p$. It is also clear from (1) that the constraints do not have any impact on the system dynamics. However, as mentioned in the Introduction, being able to converge to a desired SEP, $x_s$, after a disturbance is not sufficient; the trajectory must also *not enter* the infeasible region $\{x|h_k(x,p) \leq 0 \; \exists k \in [1,n_h]\}$. Therefore, the constraints *do* play an important role in defining the set of all desirable (stable + feasible) trajectories. This set will be referred to as the constrained stability region (CSR) of $x_s$, denoted by $A_c(x_s)$, while the corresponding constrained stability boundary will be denoted by $\partial A_c(x_s)$. Loparo [14] characterized the constrained stability boundary of DAE-governed systems with inequality constraints. In this section, we will describe the stability boundary for ODE systems.

#### A. Transformed Unconstrained System and Pseudo Unstable Equilibrium Points (UEPs)

The constrained system given in (1) can be transformed to an equivalent unconstrained system [14] as shown below.
$$\dot{x}_{n\times 1} = H(x,p)_{1\times 1} \times f(x,p)_{n\times 1} \qquad (2)$$

The unconstrained system given in (2) has the same stability region (SR) and associated boundary as the CSR of the original constrained system given in (1) but is easier to analyze due to its unconstrained nature. Therefore, we use the system described by (2) to understand the *qualitative nature of the constrained stability boundary*. In (2), $H(x,p)$ is a scalar function arrived at by multiplying all the inequality constraints. Furthermore, this system has an interplay of stability and feasibility reflected in its dynamics because of the product of $H(x,p)$ and the vector field, $f(x,p)$, in the right-hand side (RHS) of (2). When inside the feasibility region, this product only changes the length of $f(x,p)$. However, when any individual feasibility constraint is violated, i.e. $H(x,p)$ becomes negative, the product reverses the direction of $f(x,p)$. This means that the points on the feasibility boundary now also serve as equilibrium points (EPs) of this system, which we will refer to as *pseudo EPs*. The pseudo EPs will be denoted by $x_H$ to distinguish them from the original system's EPs, $x_e$. We now linearize (2) to understand the nature of pseudo EPs.

$$\Delta\dot{x} = \begin{bmatrix} \frac{\partial H}{\partial x_1}f_1(x,p) & \cdots & \frac{\partial H}{\partial x_n}f_1(x,p) \\ \vdots & \ddots & \vdots \\ \frac{\partial H}{\partial x_1}f_n(x,p) & \cdots & \frac{\partial H}{\partial x_n}f_n(x,p) \end{bmatrix} \times \Delta x \qquad (3)$$
$$+ H(x,p) \times \frac{\partial f(x,p)}{\partial x} \times \Delta x$$

The second term in RHS of (3) becomes 0 since $H(x,p) = 0$ for pseudo EPs. The connected components of $x_H$ represented by $\aleph_{x_H}$ is a $(n-1)$ dimensional manifold and thus the state matrix in the first term on RHS has $(n-1)$ eigenvalues as 0. Thus, the only possible non-zero eigenvalue equals the trace of this matrix, given by $\sum_{i=0}^{n} \frac{\partial H}{\partial x_i} f_i(x,p) = \dot{H}(x,p)$. Therefore, *a pseudo EP is stable* $(x_H^s)$, *(in the sense of Lyapunov), if $f(x,p)$ points towards the feasibility boundary, and unstable* $(x_H^u)$, *if it points away*, where the feasibility boundary serves as the local center manifold. Points on the feasibility boundary that have $\dot{H}(x,p) = 0$, do not belong to the above category, and correspond to points for which $f(x,p)$ is tangential, as shown in Figure 1. These points are referred to as semi-saddle points $(x_H^0)$, and they lie on the separating boundary between the pseudo SEP and pseudo UEP.

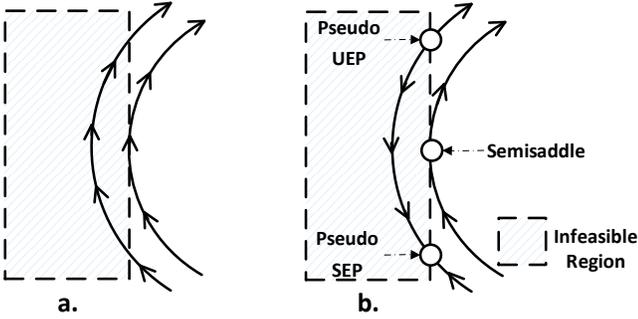

*Figure 1 Original (a) vs Transformed (b) System Dynamics*

The connected component of semi-saddle points, $\aleph_{x_H^0}$, defined by the set $\{x | H(x,p)_{1\times 1} = 0, \dot{H}(x,p)_{1\times 1} = 0, \frac{\partial H(x,p)}{\partial x} \neq 0, \ddot{H}(x,p)_{1\times 1} \neq 0\}$ has dimension of $n-2$ and was referred to as the "nice set" in [14]. There is another set of points on the feasibility boundary called the "bad set", denoted by $x_H^B$, which also has $\dot{H}(x,p)_{1\times 1} = 0$. However, at these points, either $f(x,p)$ is tangential to feasibility boundary just like semi-saddle points but $\ddot{H}(x,p)_{1\times 1} = 0$, or it is not tangential but $\frac{\partial H(x,p)}{\partial x} = 0$ indicating the presence of a local extremum point on the feasibility boundary. The connected components of bad set of points are generally of low dimension and therefore do not lie on the quasi stability boundary.

It is easy to see from (3) that each point on $\aleph_{x_H^0}$ has an $n$ dimensional center manifold (all the eigenvalues are 0). Of prime importance when characterizing the stability boundary of (2) is a subset of the center manifold $W^c(\aleph_{x_H^0})$ defined as $W^{c-}(\aleph_{x_H^0}) = \{x | H(\varphi^{(2)}(x,t,p),p) > 0 \; \forall t \geq 0, \lim_{t\to\infty} \varphi^{(2)}(x,t,p) \to \aleph_{x_H^0}\}$ where $\varphi^{(2)}$ denotes flow of system (2). Now, $\aleph_{x_H^0}$ itself is a $n-2$ dimensional sub-manifold of $R^n$ as it is given by the intersection of zero level set of two functions with all the singular points removed. We know that inside the feasible region, the system defined by (2) is merely a time scaled version of the one defined by (1) and therefore the set $W^{c-}(\aleph_{x_H^0})$ for (1) is the same as that for (2). For the system defined by (1), $W^{c-}(\aleph_{x_H^0})$ is defined as the set of points we get by back-extending in time $\aleph_{x_H^0}$ using flow of (1) and then excluding $\aleph_{x_H^0}$ i.e. $\{x | \exists (x_0, t) [(x_0, t) \in \aleph_{x_H^0} \times R_-^1 \wedge \varphi^{(1)}(x_0, t, p) = x \wedge H(\varphi^{(1)}(x_0, t_1, p)) > 0 \; \forall t_1 \in [t, 0)]\}$ where $\varphi^{(1)}$ denotes flow of system (1) and $R_-^1$ represents the negative open half space in $R^1$ i.e. $\{y \in R^1 | y < 0\}$. Since no other critical points lie on $W^{c-}(\aleph_{x_H^0})$, each point on $W^{c-}(\aleph_{x_H^0})$ can be defined by a unique combination of starting point on $\aleph_{x_H^0}$ and the time $t$ it takes for the reverse flow of (1) to reach that point. Also, for fixed $p$, $\varphi^{(1)}: \aleph_{x_H^0} \times t \to W^{c-}(\aleph_{x_H^0})$ is a diffeomorphism. Here, $t \in R_-^1$ which is diffeomorphic to $R^1$ [15]. Furthermore, since $\aleph_{x_H^0}$ is a zero level set of smooth functions $H(x,p)$ and $\dot{H}(x,p)$ without any points having $\frac{\partial H(x,p)}{\partial x} = 0$ or $\frac{\partial \dot{H}(x,p)}{\partial x} = 0$ by definition, by implicit function theorem, it is a smooth $n-2$ dimensional submanifold of $R^n$. Therefore, it is straightforward to say that $W^{c-}(\aleph_{x_H^0})$ is an $n-1$ dimensional smooth manifold.

### B. Characterization of Quasi-Stability Boundary

The following assumptions must be satisfied for the stability boundary characterization of (2) [14]:

(A1) No $x_e$ on the feasibility boundary.

(A2) $\aleph_{x_H^0}$ is *topologically* dense in the set $\{x | H(x,p) = 0, \dot{H}(x,p) = 0\}$ and $\aleph_{x_H^B}$ has a maximum dimension of $n-3$.

(A3) All original system EPs and periodic orbits on the stability boundary are hyperbolic.

(A4) $W^s(x_e)$, $W^s(\aleph_{x_H^B})$, and $W^{c-}(\aleph_{x_H^0})$ intersect transversally with $W^u(x_e), W^u(\aleph_{x_H^B})$ and $W^{c+}(\aleph_{x_H^0})$. However, $W^{c-}(\aleph_{x_H^0})$ and $W^{c+}(\aleph_{x_H^0})$ are not transversal for the same $\aleph_{x_H^0}$.

(A5) Any trajectory on the stability boundary converges to one of the EPs or periodic orbits on the boundary.

The stability boundary structure of a generalized nonlinear system can be very complex. On the contrary, the *quasi-stability region is a practical SR whose boundary is the boundary of closure of $A(x_s)$, denoted by $\partial \bar{A}(x_s)$* [16]. Taking the closure removes the low dimensional components of the stability boundary, which considerably simplifies the analysis. It has been shown in [14] that for constrained systems satisfying the assumptions (A1)-(A5), the quasi stability boundary of (2) comprises stable manifolds of type 1 $x_e^u$, type 2 periodic orbits, $W^{c-}$ of $(n-2)$ dimensional $\aleph_{x_H^0}$ whose $W^{c+}$ manifolds intersect the SR, and unstable portions of the feasibility boundary $\aleph_{x_H^u}$. To better visualize the stability boundary of a system of the form (2), we plot the SR of a SMIB system with inequality constraint, $H(\delta, \omega) = 2 - \delta > 0$, as shown in Figure 2. The feasibility boundary is shown by red and blue dotted line with the blue component representing $\aleph_{x_H^s}$ and red component representing $\aleph_{x_H^u}$. $\aleph_{x_H^0}$ is shown as a green dot between these two components with $W^{c-}(\aleph_{x_H^0})$ serving as a part of the stability boundary which is the boundary of the whole colored region plotted. There is also a type 1 $x_e^u$ along with its stable manifold serving as a part of the stability boundary on the left.

It is important to point out here that the transformed system defined by (2) was only used to familiarize oneself with how constraints influence the characteristics of the CSR. The CCT sensitivity derivations described in the next section will be performed on the original (untransformed) constrained system defined by (1).

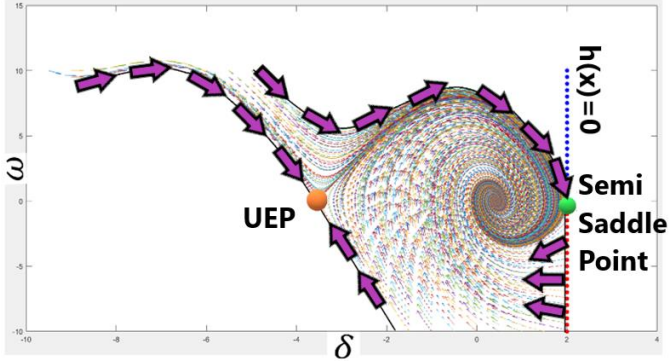

Figure 2 Example Constrained System's Stability Region

## IV. SENSITIVITY DERIVATIONS

In the subsequent sections, superscripts will be used to define the network topology status, namely, pre-fault, fault-on, and post-fault. For e.g., $f^{fault}$ will represent the vector field of fault-on system. A critical fault-on trajectory for a given fault would be one that intersects the constrained stability boundary of the post-fault system i.e. $\partial A_c(x_s^{post})$. There are two ways how this could happen: (i) it intersects the $W^s(x_e^{post})$ or $W^{c-}(\aleph_{x_H^0}^{post})$ or, (ii) it intersects the feasibility boundary $\aleph_{x_H^s}^{fault}$ or $\aleph_{x_H^u}^{post}$. The two components of the second category can be combined by multiplication: $\{x|H^{fault}(x,p) \times H^{post}(x,p) = 0\}$. Note that unlike the previous section, here we will analyze the effect of variation of parameter $p$ on the CCT sensitivity.

The parametric constrained stability boundary for a fixed $p$ can be written in the form $S(x,p) = 0$ due to its $n-1$ dimensional nature ($n$ dimensional in $x-p$ space). The intersection of parametric fault trajectory and constrained stability boundary exists under parameter changes if they intersect transversally [17]. Thus, for the same change of $p$, the CCT, $t_{cr}$, would be changed such that the new state at the fault clearing time lies on the new constrained stability boundary. In order to achieve this, we need to derive the sensitivity formula for the state value at the beginning of fault on trajectory denoted by $x_0$, at the time of fault clearing, $t_{cl}$, denoted by $\varphi^{fault}(x_0, t_{cl}, p)$, and that of the associated relevant portion [18] of the constrained stability boundary which can be represented as a zero level set of a function as $S(x,p) = 0$ due to it being $n-1$ dimensional. This as illustrated in Figure 3. The original fault trajectory is shown by solid blue line and the perturbed version is shown by dashed blue line. The double-headed arrows connect points along both the trajectories at the same time. The perturbed and original constrained stability boundary are shown by black solid and black dashed lines, respectively. In this example, on perturbing $p$, the constrained stability boundary moves further away resulting in an increase of $t_{cl}$, and therefore, $t_{cr}$. This shows a positive sensitivity of $t_{cr}$ to $p$.

As discussed in Section III. , there are three structurally distinct portions of the $\partial \overline{A_c}$, namely, unstable portions of the feasibility boundary, manifold defined by trajectories converging to an $n-2$ dimensional connected component of semi-saddle points on the feasibility boundary, and stable manifold of type 1 UEPs and type 2 limit cycles. As limit cycles rarely occur in a power system, they are not considered in this analysis. Depending on the mode of loss of stability/feasibility of a given critical fault trajectory, the sensitivity of the appropriate relevant portion of the boundary must be calculated. We will now present the derivations for each case. In the following sub-sections, $p$ is assumed to be scalar. The sensitivity is computed for the base critical trajectory (obtained from time domain simulations) having parameter value, $p_0$. For sake of clarity, some of the constants coefficients in the equations will be replaced by new variables.

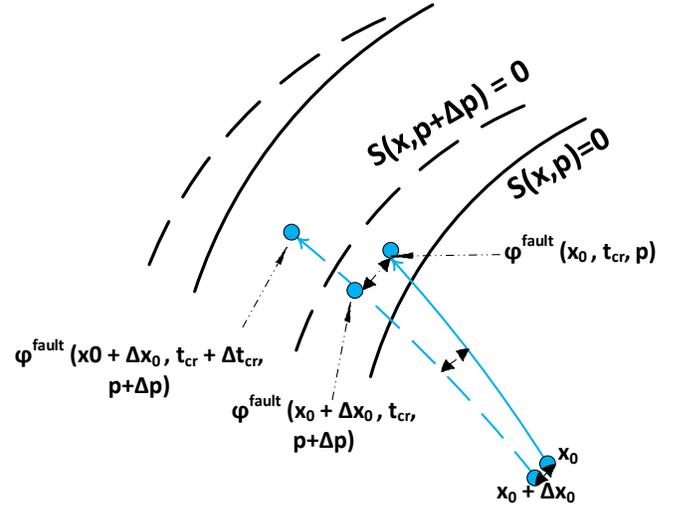

Figure 3 Overall Process

### A. Sensitivity of the State Value at Fault Clearing ($x_{cl}$)

Usually, the starting point of the parametric fault-on trajectory is the SEP of the pre fault system, $x_s^{pre}(p)$, which lies on a single dimensional manifold given by $f^{pre}(x,p) = 0$. Let the starting point of base critical trajectory, $x_s^{pre}(p_0)$, be denoted by $x_{s_0}^{pre}$. Sensitivity of $x_s^{pre}(p)$ evaluated at $x_{s_0}^{pre}$ is,

$$\frac{\Delta x_s^{pre}(p)}{\Delta p}\bigg|_{p_0} = M_{4(n \times 1)} = -\left[\frac{\partial f^{pre}(x,p)}{\partial x}\bigg|_{x_{s_0}^{pre}, p_0}\right]^{-1} \times \frac{\partial f^{pre}(x,p)}{\partial p}\bigg|_{x_{s_0}^{pre}, p_0} \quad (4)$$

Next, we compute the sensitivity of the state value at any general fault clearing time, $t_{cl}$, for a fault trajectory starting from $x_s^{pre}(p)$. This point is a function of $t_{cl}$, $p$, and $x_s^{pre}$ and is therefore denoted by $x_{cl}(x_s^{pre}(p), t_{cl}, p)$. Since $x_s^{pre}$ itself is a function of $p$, $x_{cl}$ effectively lies on a two-dimensional manifold in the $x-p$ space. Let the CCT of the base critical fault trajectory be $t_{cr}$ and the state value at that time be $x_{cr}$. Calculating sensitivity of $x_{cl}$ around the base critical trajectory at $t_{cr}$ and combining with (4), we get,

$$\frac{\Delta x_{cl}(x_s^{pre}(p), t_{cl}, p)}{\Delta p}\bigg|_{x_{s_0}^{pre}, t_{cr}, p_0} \quad (5)$$

$$= M_1 \times M_4 + M_2 \times \frac{\Delta t_{cl}}{\Delta p}\bigg|_{p_0} + M_3$$

where, $M_{1(n \times n)} = \frac{\partial \varphi^{fault}(x_0,t,p)}{\partial x_0}\Big|_{x_{s_0}^{pre},t_{cr},p_0}$, $M_{2(n \times 1)} = \frac{\partial \varphi^{fault}(x_0,t,p)}{\partial t}\Big|_{x_{s_0}^{pre},t_{cr},p_0} = f^{fault}(x,p)\Big|_{x_{cr},p_0}$, and $M_{3(n \times 1)} = \frac{\partial \varphi^{fault}(x_0,t,p)}{\partial p}\Big|_{x_{s_0}^{pre},t_{cr},p_0}$.

### B. Sensitivity of Combined Feasibility Boundary of Fault-on and Post-Fault System

The direct loss of feasibility (not requiring integration of post-fault trajectory) for constrained power systems happens if the sustained fault trajectory intersects either the (i) stable component of the feasibility boundary of the fault-on system w.r.t. $f^{fault}(x,p)$, denoted by $\aleph_{x_H^s}^{fault}$ or (ii) feasibility boundary of the post-fault system which is unstable w.r.t. $f^{post}(x,p)$, denoted by $\aleph_{x_H^u}^{post}$. These two components can be combined together by multiplication: $\{x|H^{comb}(x) = H^{fault}(x) \times H^{post}(x) = 0\}$, to get the combined intersecting boundary. There should not be constraint functions present in both $H^{fault}(x)$ and $H^{post}(x)$ as it may make $H^{comb}(x)$ positive definite in some regions. This usually arises in situations when the constraints are independent of network topology, so same constraints occur in all network topologies. In cases where $x_{cr}$ satisfies $H^{comb}(x_{cr},p_0) = 0$, for small variation in $p$, $t_{cl}$ should change such that $H^{comb}(x_{cl}(x_s^{pre}(p),t_{cl},p),p)$ remains 0 in order for $t_{cl}$ to still represent CCT under parameter changes. This can be written as,

$$\frac{\Delta x_{cl}}{\Delta p}\Big|_{x_{s_0}^{pre},t_{cr},p_0} = M_5^- \times M_6 \tag{6}$$

where, $M_{5(1 \times n)} = \frac{\partial H^{comb}(x,p)}{\partial x}\Big|_{x_{cr},p_0}$, and $M_{6(1 \times 1)} = -\frac{\partial H^{comb}(x,p)}{\partial p}\Big|_{x_{cr},p_0}$

Substituting (5) in (6) yields the change in $t_{cl}$ required w.r.t change in $p$, which is also the CCT sensitivity for the given mode of loss of stability/feasibility,

$$\boxed{\frac{\Delta t_{cl}}{\Delta p}\Big|_{p_0} = [M_5 \times M_2]^{-1} \times (M_6 - (M_5 \times (M_1 \times M_4 + M_3)))} \tag{7}$$

### C. Sensitivity of Post-Fault Trajectory's End Point

CCT sensitivities in the remaining categories, namely, intersecting feasibility boundary in the post-fault phase or intersecting stable manifold of controlling UEP, are derived by combining sensitivity of the end point of the post-fault trajectory with the sensitivity of the local characterization of the constrained stability boundary around that point. Now, the constrained stability boundary also contains the stable manifold, $W^s$, of some original system's EPs and $W^{c-}$ of connected component of semi-saddle pseudo EPs. $W^s$ and $W^{c-}$ can be visualized as a surface of adjacent trajectories reaching the same set of points [19]. Since a critical fault trajectory is a single dimensional manifold, it would be intersecting the stability boundary at a single point. Therefore, we will focus on the emerging post-fault trajectory from that point, referred to as the *critical post-fault trajectory*. The sensitivity will be derived around the base critical post-fault trajectory.

The local characterization of the constrained stability boundary at the end point of the post-fault trajectory is normally available/can be derived. Therefore, the first step is to compute the sensitivity of the post-fault trajectory's end point which is a function of the starting point of the post-fault trajectory/ending point of the fault-on trajectory, i.e. $x_{cl}$ for ODE systems, time elapsed along the post-fault trajectory $t_{end}$, and $p$; therefore, it can be written as $x_{end}^{post}(x_{cl},t_{end},p)$. Here, $t_{end}$ is assumed to be large enough such that $x_{end}^{post}$ reaches the region in state space where the known local characterization of the constrained stability boundary, denoted by $S^{loc}(x,p) = 0$, holds true. Let $T$ denote the time it takes for the base critical trajectory to come arbitrarily close to a point at which local characterization of the constrained stability boundary is available. Also, let the value of state variable at that time be $x_T$ i.e. $x_T = \varphi^{post}(x_{cr},T,p_0)$. First, evaluating the sensitivity of $x_{end}^{post}$ at that point, we get,

$$\frac{\Delta x_{end}^{post}(x_{cl},t_{end},p)}{\Delta p}\Big|_{x_{cr},T,p_0} \tag{8}$$

$$= O_1 \times \frac{\Delta x_{cl}}{\Delta p}\Big|_{x_{s_0}^{pre},t_{cr},p_0} + O_2 \times \frac{\Delta t_{end}}{\Delta p}\Big|_{p_0} + O_3$$

where, $O_{1(n \times n)} = \frac{\partial \varphi^{post}(x_0,t,p)}{\partial x_0}\Big|_{x_{cr},T,p_0}$, $O_{2(n \times 1)} = \frac{\partial \varphi^{post}(x_0,t,p)}{\partial t}\Big|_{x_{cr},T,p_0} = f^{post}(x,p)\Big|_{x_T,p_0}$, and $O_{3(n \times 1)} = \frac{\partial \varphi^{post}(x_0,t,p)}{\partial p}\Big|_{x_{cr},T,p_0}$.

Substituting (5) in (8), we get,

$$\frac{\Delta x_{end}^{post}(x_{cl},t,p)}{\Delta p}\Big|_{x_{cr},T,p_0} \tag{9}$$

$$= [O_1 \times M_2 \quad O_2] \times \begin{bmatrix} \frac{\Delta t_{cl}}{\Delta p}\Big|_{p_0} \\ \frac{\Delta t_{end}}{\Delta p}\Big|_{p_0} \end{bmatrix} + O_1 \times (M_1 \times M_4 + M_3) + O_3$$

In the following sub-sections, we will derive the sensitivity of different local characterizations of constrained stability boundary (depending on mode of loss of stability/feasibility) on which $x_{end}^{post}$ should stay for the overall trajectory to remain critical.

### D. Sensitivity of Attracting Subset of Center Manifold of $n-2$ Dimensional Connected Component of Semi-Saddle Pseudo EP

The CCT sensitivity formula derived in this section will be used if the base critical trajectory after some time, $T$, along the post-fault system defined by (1) grazes the feasibility boundary i.e. $H^{post}(x_T,p_0) = \frac{\partial H^{post}(x_T,p_0)}{\partial x} \times f^{post}(x_T,p_0) = 0$. It is important to note here that for the transformed system given by (2), $T$ could $\to \infty$ but for (1), it is finite due to the absence of any critical points in the form of feasibility boundary. As discussed in Section III. A. , after removing the points belonging to the "bad set", the local characterization of constrained stability boundary at any point on the connected component of semi-saddle pseudo EPs of the parametric post-fault system is given by $\{H^{post}(x,p) = 0, \frac{\partial H^{post}(x,p)}{\partial x} \times f^{post}(x,p) = 0\}$. Therefore, $x_{end}^{post}$ must satisfy this constraint for the post-fault trajectory to remain critical, under variations in $p$, as demonstrated in Figure 4. Evaluating the sensitivity at end point ($x_T$) of the base critical post-fault trajectory, we get,

$$O_4 \times \frac{\Delta x_{end}^{post}(x_{cl},t,p)}{\Delta p}\bigg|_{x_{cr},T,p_0} = O_5 \quad (10)$$

where, $O_{4_{(2\times n)}} = \begin{bmatrix} \frac{\partial H^{post}(x,p)}{\partial x}\big|_{x_T,p_0} \\ \frac{\partial \left[\frac{\partial H^{post}(x,p)}{\partial x} \times f^{post}(x,p)\right]}{\partial x}\bigg|_{x_T,p_0} \end{bmatrix}$ and

$$O_{5_{(2\times 1)}} = -\begin{bmatrix} \frac{\partial H^{post}(x,p)}{\partial p}\big|_{x_T,p_0} \\ \frac{\partial \left[\frac{\partial H^{post}(x,p)}{\partial x} \times f^{post}(x,p)\right]}{\partial p}\bigg|_{x_T,p_0} \end{bmatrix}$$

Combining (9) and (10), we get the expression for CCT sensitivity as,

$$\boxed{\begin{bmatrix} \frac{\Delta t_{cl}}{\Delta p}\big|_{p_0} \\ \frac{\Delta t_{end}}{\Delta p}\big|_{p_0} \end{bmatrix} = [O_4 \times [O_1 \times M_2 \quad O_2]]^{-1} \times (O_5 \\ - O_4 \times (O_1 \times (M_1 \times M_4 + M_3) \\ + O_3))} \quad (11)$$

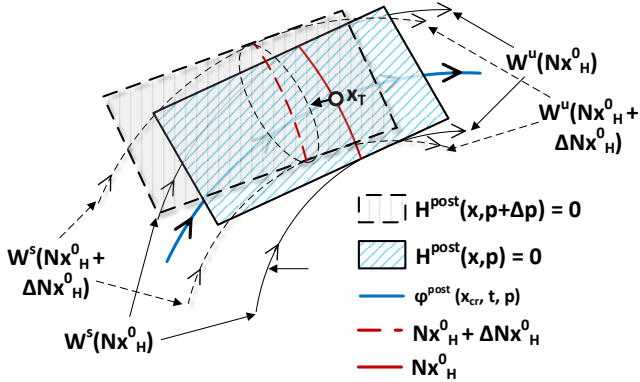

*Figure 4 Semi-Saddle Surface Sensitivity*

### E. Sensitivity of Stable Manifold of Type 1 CUEP

The CCT sensitivity formula derived in this section will be used if the base critical post-fault trajectory does not return to the SEP, which is usually marked by loss of synchronism or voltage collapse. Since the inequality constraints do not impact the dynamics, they do not come into play, and hence can be ignored. To be specific, the sustained fault trajectory exits the SR of the original system without any constraints through the stable manifold of a UEP which is referred to as the controlling UEP (CUEP) [16] that locally defines the constrained stability boundary. For constrained quasi-stability boundary, we are only interested in type 1 CUEPs whose $n-1$ dimensional stable manifold makes up the constrained stability boundary. Constant energy surface passing through the CUEP has been traditionally used as an approximation to the stability boundary with its sensitivity used to estimate the sensitivity of the CCT. Instead, here we will use the local characterization of type 1 CUEP's stable manifold, $W^s$, along with the sensitivity of the post-fault trajectory's end point presented in Section IV. C.

In this section, $T$ is chosen appropriately for the base critical trajectory such that $x_T$ is in close proximity to $x_{cu_0}^{post} = x_{cu}^{post}(p)\big|_{p_0}$, say inside some open volume $B$ containing $x_{cu_0}^{post}$. For a type $k$ UEP, the stable manifold of an EP is locally characterized by a hyperplane spanned by $n-k$ stable eigenvectors [19]. The local characterization of the stable manifold of type-1 $x_{cu}^{post}(p)$ is available at the $x_{cu}^{post}(p)$ and is defined by the hyperplane spanned by the stable eigenvectors. However, $T$ needs to tend to $\infty$ such that $x_T \to x_{cu}^{post}(p)$ for that characterization to be used. An effective approach was proposed in [13] to get a local characterization of $W^s(x_{cu}^{post}(p))$ at $x_T$; that approach will be used here directly. More details about the approach can be found in [13].

The volume $B$ mentioned previously which contained $x_{cu_0}^{post}$ is defined such that a differentiable parameter varying chart $\gamma: B \times p \to R^n$ exists such that the parametric stable manifold $W^s(x_{cu}^{post}(p))$ is given by [15],

$$S^{loc}(x,p) = w(p) \times (\gamma(x,p) - \gamma(x_{cu}^{post}(p),p)) = 0 \quad (12)$$

where $w(p)$ is the left eigen vector corresponding to the only unstable eigenvalue of $x_{cu}^{post}(p)$. Now, flow at non-critical points is a diffeomorphism, so we extend the definition of $S^{loc}(x,p)$ which was defined around $x_T$ inside $B$ to other points around the base case critical post-fault trajectory under small variations in $p$ as follows,

$$S^{loc}(x,p) = w(p) \times (\gamma(\varphi^{post}(x,\tau(x,p),p),p) \\ - \gamma(\varphi^{post}(x_{cl},T,p),p)) = 0 \quad (13)$$

where $\tau(x,p)$ denotes the time taken for a parametric post-fault trajectory starting at an arbitrary $x$ to reach $B$. Now, $\frac{\partial S^{loc}}{\partial x} \times f^{post}(x,p) = 0$ since $\{x|S^{loc}(x,p) = 0\}$ is an invariant set. Next, we evaluate the sensitivity of $x_{cu}^{post}(p)$ to $p$ at $x_{cu_0}^{post}, p_0$. As $x_{cu}^{post}(p)$ lies on a single dimensional manifold satisfying $f^{post}(x_{cu}^{post}(p),p) = 0$, differentiating it gives,

$$\frac{\Delta x_{cu}^{post}}{\Delta p}\bigg|_{p_0} = O_{6_{(n\times 1)}} \quad (14)$$

$$= -\left[\frac{\partial f^{post}(x,p)}{\partial x}\bigg|_{x_{cu_0}^{post},p_0}\right]^{-1} \\ \times \frac{\partial f^{post}(x,p)}{\partial p}\bigg|_{x_{cu_0}^{post},p_0}$$

Now, it has been shown in [13] that for sufficiently large values of $T$,

$$\frac{\partial S^{loc}}{\partial x}\big|_{x_T,p_0} \to w(p)\big|_{p_0}, \frac{\partial S^{loc}}{\partial p}\big|_{x_T,p_0} \to -w(p)\big|_{p_0} \times \\ \frac{\Delta x_{cu}^{post}(p)}{\Delta p}\big|_{p_0}, \frac{\partial \gamma(x)}{\partial x} = Identity \quad (15)$$

Finally, differentiating (13) and evaluating at $x_T$ for a sufficiently large value of $T$ using (15), we get,

$$w(p)\big|_{p_0} \times \left(\frac{\Delta x_{end}^{post}(x_{cl},t_{end},p)}{\Delta p}\bigg|_{x_{cr},T,p_0}\right) \\ = w(p)\big|_{p_0} \times O_6 \quad (16)$$

Substituting (9) in (16), we get,

$$\left.\frac{\Delta t_{cl}}{\Delta p}\right|_{p_0} = \frac{w(p)|_{p_0} \times (O_6 - O_3 - O_1 \times (M_1 \times M_4 + M_3))}{w(p)|_{p_0} \times O_1 \times M_2} \quad (17)$$

## V. Extension to Constrained Differential Algebraic Equation (DAE)-Governed Systems

Parametric power system governed by DAEs with inequality constraints can be written as

$$\dot{x} = f(x, y, p) \quad (18)$$
$$0 = g(x, y, p)$$
$$h(x, y, p) \geq 0$$

The equality constraint $g = 0$ which is given for each configuration i.e. pre-fault, fault-on and post-fault, gives the corresponding surface in the overall state space on which the system evolves. The system jumps between these surfaces whenever switching happens, with the assumption that $x$ stays the same during switching while $y$ jumps. On a given surface, as long as $\frac{\partial g(x,y,p)}{\partial y}$ is invertible, $y$ can be locally written as a function of $x$ (implicit function theorem). The scenarios under which this model can break down (i.e. $\frac{\partial g(x,y,p)}{\partial y}$ becomes singular) will be explored in a future work.

As such, constrained DAE systems can be analyzed by using an extended state space representation where the value of $y$ on each configuration/surface is defined as a new state variable. However, at a given time, only the configurations/surfaces the system will switch to must be considered. For example, dynamics in the fault-on configuration in the extended representation can be written as,

$$\dot{x} = f^{fault}(x, y^{fault}, p) \quad (19)$$
$$0 = g^{fault}(x, y^{fault}, p)$$
$$h^{fault}(x, y^{fault}, p) \geq 0$$
$$0 = g^{post}(x, y^{post}, p)$$
$$h^{post}(x, y^{post}, p) \geq 0$$

Here, the evolution of $x$ is only a function of the $x$ and $y$ values on the fault surface, i.e. $y^{fault}$. Naturally, no other surfaces are analyzed once the system is in post-fault configuration since no further switching happens from it. Assuming the system starts from the pre-fault system's SEP, the starting point lies on the following manifold,

$$0 = f^{pre}(x, y^{pre}, p) \quad (20)$$
$$0 = g^{pre}(x, y^{pre}, p)$$
$$0 = g^{fault}(x, y^{fault}, p)$$
$$0 = g^{post}(x, y^{post}, p)$$

The trajectory sensitivities are straightforward to compute keeping in mind that they must be computed for the extended version of the system for a given configuration. Further, only variations in trajectories w.r.t. starting values of $x$ and $p$ must be computed, as they are the only independent quantities. Lastly, derivations that require rate of change of trajectories, $\dot{y}$, can be derived by differentiating the appropriate equality constraint linking $x$ and $y$. For example, $\dot{y}^{post}$ along fault trajectory can be obtained by differentiating the 4th term in (20). The other sensitivities can be derived in a manner similar to what was done for constrained ODE systems in Section IV.

## VI. Overall Computation

This section discusses the different computations involved in finding the sensitivity of CCT of a given fault to various parameter changes. The steps involved are as follows:

1. Computationally tractable direct method for computing CCT for constrained systems is a challenge due to changes in nature of the constrained stability boundary as compared to the traditional unconstrained formulation of power systems [5]. Hence, for the constrained system under study, CCT as well as critical fault-on and post-fault trajectories are found using time domain simulation (TDS) for fixed values of parameters, $p_0$, using the following algorithm.

    **Algorithm: CCT and Critical Trajectory Computation using TDS for Constrained Systems**
    i. **INITIALIZE** stable clearing time $t_{stable}$ and unstable clearing time $t_{unstable}$. $t_{unstable}$ is set to the time at which the sustained fault trajectory intersects the feasibility boundary.
    ii. **SET** $t_{cl} = \frac{t_{stable}+t_{unstable}}{2}$. $\varphi^{fault}(x_s^{pre}(p_0), t_{cl}, p_0)$ is denoted by $x_{cl}$.
    iii. **INITIALIZE** $t_1 = t_2 = \infty$. $x_T, T = null$
    iv. Integrate the post-fault trajectory for a long time $T_{max}$.
    v. **UPDATE** $t_1$ equal to time at which $H^{post}(x)$ crosses 0 or $H^{post}(x) \leq 1e - 5$.
    vi. **IF** $t_1 \leq T_{max}$, $t_{unstable} = t_{cl}$.
    vii. **IF** $\varphi^{post}(x_{cl}, T_{max}, p_0) \neq x_s^{post}(p_0)$, $t_{unstable} = t_{cl}$. Update $t_2$ to time where $||f^{post}(x, p_0)||$ acquires a local minimum value along the post-fault trajectory and is $\leq 1e - 3$.
    viii. **IF** $\min(t_1, t_2) < \infty, T = \min(t_1, t_2), x_T = \varphi^{post}(x_{cl}, T, p_0)$.
    ix. **IF** $|t_{stable} - t_{unstable}| \geq 0.01$ **OR** $T = null$, **GOTO** ii.
    x. **STOP**

    The following things must also be noted:
    a. The transformed unconstrained system given in (2) can also be used for TDS. However, an adaptive step size will be required for simulation as the time scale varies drastically with the value of $H(x, p)$ along a trajectory. This would require using a solver for stiff systems, which would then increase the TDS computation.
    b. It is very difficult to find the exact time at which a fault trajectory intersects $W^s(CUEP)$. Therefore, we use the approach used for finding CUEP for gradient system in the boundary controlling unstable (BCU) method [16].

2. Besides the various Jacobian computations, the following sensitivities must be computed as part of the overall process using [1]:
    i. Integrating the fault on trajectory till $t_{cr}$ to compute $\left.\frac{\partial \varphi^{fault}(x_0, t, p)}{\partial x_0}\right|_{x_{s_0}^{pre}, t_{cr}, p_0}, \left.\frac{\partial \varphi^{fault}(x_0, t, p)}{\partial p}\right|_{x_{s_0}^{pre}, t_{cr}, p_0}$
    ii. If the loss of stability/feasibility is not intersection of fault-on trajectory with the feasibility boundary, compute, $\left.\frac{\partial \varphi^{post}(x_0, t, p)}{\partial x_0}\right|_{x_{cr}, T, p_0}, \left.\frac{\partial \varphi^{post}(x_0, t, p)}{\partial p}\right|_{x_{cr}, T, p_0}$.

When using the proposed approach on large-scale systems, the main bottleneck is the computation of trajectory sensitivities. This has been overcome by using parallel programming and sparsity techniques as done in [20], [21]. Furthermore, the characterization of stability boundary on which our derivations are based hold true in general regardless of the size of the system.

## VII. RESULTS

We will use the following categories to denote the instability/infeasibility phenomenon:

**1**: fault trajectory directly intersects the feasibility boundary
**2**: post-fault trajectory intersects the feasibility boundary
**3**: post-fault trajectory does not return to $x_s^{post}$

### A. Single-Machine Infinite-Bus (SMIB) System Results

It must also be noted here that although a SMIB system is analyzed here for ease of understanding, the proposed methodology is applicable to larger systems as well. The dynamics of an SMIB system is described by the state equation,

$$\dot{\delta} = \omega \quad (21)$$
$$M\dot{\omega} = P_m - \frac{EV}{X}\sin(\delta) - D \times \omega$$

Here, $\delta$ denotes rotor angle, $\omega$ denotes angular speed deviation, $M$ is inertia, $P_m$ is mechanical power input, $D$ is damping, $E$ is internal emf of the generator, $V$ is voltage of the infinite bus, and $X$ is the total impedance. Fault being analyzed is on the infinite bus, i.e. $\frac{EV(fault)}{X} = 0$ and cleared without changing the topology. The constraints assumed are of the form $h(x) = [\delta^{max} - \delta, \omega^{max} - \omega]^T$ arising from out-of-step relay setting for the generator, and frequency threshold from over-frequency ride through limit on some large RG in that area. The fixed parameter values are $D = 0.5$, $\frac{EV(pre)}{X} = \frac{EV(post)}{X} = 1$. Sensitivities are computed at various parameter value combinations, where $p = [P_m, M, \delta_{max}, \omega_{max}]$.

We first analyze the effect of generator mechanical input $P_m$ on CCT at a given operating point. This is important to study as it represents the change in dispatch. The mode of loss of stability/feasibility at the point of study ($p = [0.6, 0.25, 2.4434, 1]$) is Category **1** as the sustained fault trajectory directly intersects the feasibility boundary as $P_m$ varies. Figure 5 shows the actual CCT vs $P_m$ obtained through TDS. The red dotted lines depict the computed CCT sensitivities using the sensitivity formula derived in Section IV. B. at different initial parameter values. It can be seen from the figure that the dotted lines are tangential to the original curve, which proves the validity of the formula.

Next, we try to understand the implications of changing inertia on meeting frequency constraint $\omega^{max}$. It can be seen from Figure 6 that as the inertia increases, the fault needs to be sustained longer to violate the frequency limits, a phenomenon that is expected. The trend stays the same up to a certain extent but then suddenly changes due to a change in the instability phenomenon. In this case, for $M \in [0.1: 0.15]$, the sensitivity is calculated using the derivation in Section IV. B. (Category **1**), while for $M \in (0.15, 0.3]$, it is computed using Section IV. D. (Category **2**). The sensitivity estimates continue to be tangential, thereby validating the formulae.

To analyze this phenomenon further, the relevant portion of CSRs for the post-fault constrained system are plotted using black arrows under inertia variation in Figure 7. The thick black line starting from green point (pre-fault system's SEP) to orange point (exit point) represents the sustained fault trajectory in each case. The two portions of the feasibility boundary defined by $\delta = \delta^{max}$ and $\omega = \omega^{max}$ are shown by vertical and horizontal dotted lines, respectively. Inertia plays the role of reducing the effect of angular excursion on speed, which can be seen from the changing shape of CSRs. For higher inertia values, at same speed, larger angle deviation/synchronization torque is needed to stabilize the system. For the given constrained system, we can see that for low inertia values of [0.1,0.15], the top portion of the feasibility boundary makes up the relevant portion of the constrained stability boundary since $\omega$ excursions are higher for the same fault and therefore more likely to be violated. As the inertia is increased above 0.15, the sustained fault trajectory's direction becomes more horizontal since inertia does not let the speed grow quickly. Therefore, the fault trajectory now switches to violating the angle constraint in the post-fault phase, which then becomes the new mode of infeasibility. That is, the relevant portion of the constrained stability boundary is now the stable manifold of the semi-saddle pseudo EP (blue ball) on the feasibility boundary portion given by $\delta = \delta^{max}$. Thus, it is the structural change in the relevant portion of the constrained stability boundary that causes a sharp change in the CCT vs $M$ plot shown in Figure 6.

Finally, we compute the sensitivities to the angle constraint $\delta_{max}$ as shown in Figure 8. They are also found to be tangential, thereby, confirming our derivation. For this case also, the constrained stability boundary undergoes a structural change but of a different type as compared to the last case. A closer look at the CSR in Figure 9 shows that the mode of loss of stability/feasibility for $\delta_{max}$ values of 1.74–2.09 is $\delta$ crossing the portion of the feasibility boundary defined by $\delta = \delta_{max}$ in the post-fault phase (Category **2**). The relevant semi-saddle pseudo EP is denoted by blue circle. The CCT increases with $\delta^{max}$ since longer fault clearing times result in larger angular excursions and therefore increasing chance of violating $\delta^{max}$ limit. The given situation would normally be seen for a conservatively set out-of-step relay setting where the relay trips without the generator actually going out-of-step. However, as this constraint is relaxed further, one of the original system's UEP (yellow circle) crosses the $\delta = \delta^{max}$ portion of the feasibility boundary and enters the feasible region, a phenomenon that is similar to the *singularity-induced bifurcation*. The stable manifold of that UEP now starts to serve as the relevant portion of the constrained stability boundary for the fault under study. Therefore, the instability mechanism now becomes loss of synchronism (Category **3**) with the $\delta^{max}$ limit being violated after the synchronism is already lost. This is usually how out-of-step relays are set when they are made to wait for the angles to grow large enough before acting. On further increasing $\delta_{max}$ beyond 2.1, the CCT does not change at all since that portion of the feasibility boundary no longer matters.

It must also be pointed out that in conventional unconstrained power systems the parameters under study usually affect the overall system dynamics, including fault trajectory and all portions of the post-fault system constrained stability boundary. This makes the relevant portion of the constrained stability boundary structurally stable and consequently the CUEP varies smoothly with parameter changes. For constrained systems, the parameters that define the constraints only affect one or more portions of the feasibility boundary and not the system dynamics itself. This means that with those parameters, all components of the constrained stability boundary do not vary, and only portions related to the constraints vary. This makes the relevant portion of the constrained stability boundary more prone to structural changes as seen in the two previous scenarios.

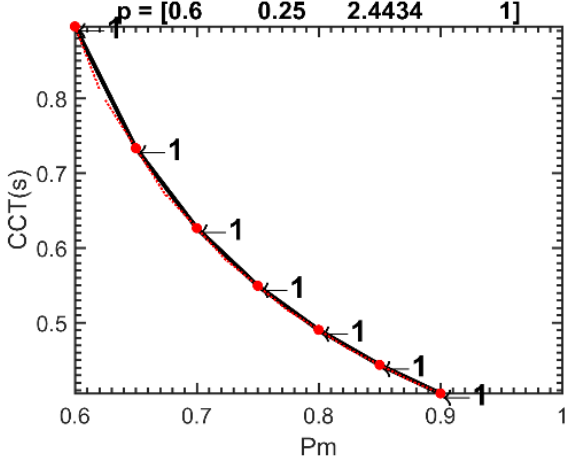

Figure 5 Single Machine System CCT vs $P_m$

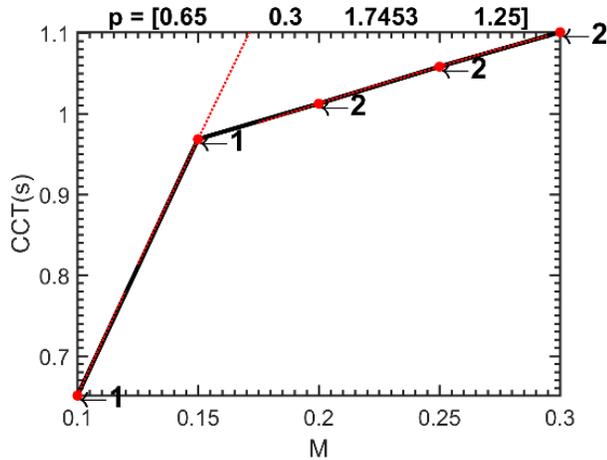

Figure 6 Single Machine System CCT vs M

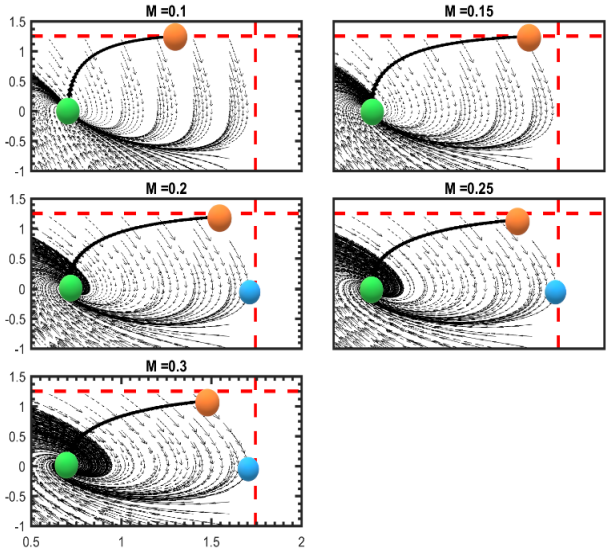

Figure 7 Single Machine System Changing CSR with M

## B. Multi-Machine System

Next, we validate the methodology on a multi-machine system [22] shown in Figure 10 with the network buses reduced to internal generator buses and assuming uniform damping coefficient $\frac{D}{M} = 4$. The fault being studied is on line 1-2 on the bus 1 side which is cleared by disconnecting the line. First, the variation of $P_m$ of generator 1 is studied. The constraint being analyzed is $\delta_1 - \delta_2 \leq \frac{\pi}{2}$. In this system, the critical trajectory eventually violates this constraint after clearing the fault (Category **2**). As $P_m$ is increased, for the same fault, generator 1 accelerates much more due to a bigger input-output imbalance (as seen in Figure 11), which results in an increased tendency to violate this angle difference constraint.

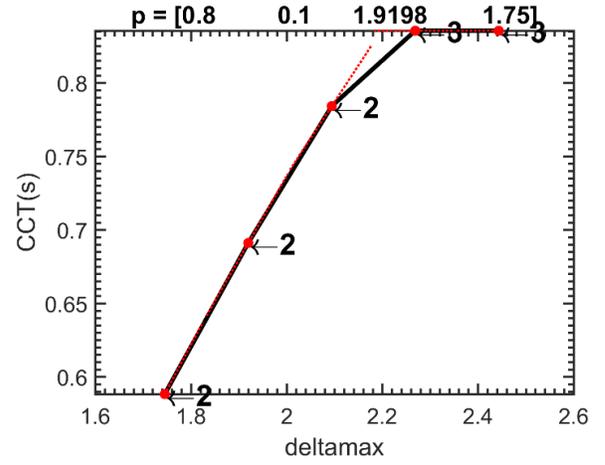

Figure 8 Single Machine System CCT vs $\delta_{max}$

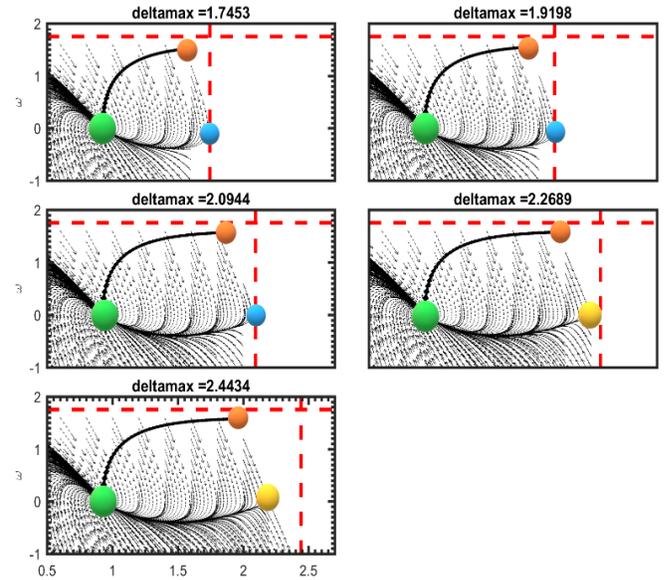

Figure 9 Single Machine System CSR for $\delta_{max} = [1.74, 2.09, 2.26, 2.44]$

Next, we study the CCT variation with increase in $P_m$ of machine 2. As seen in Figure 12, the CCT gradually increases in this case indicating a decreased tendency to violate the angle difference constraint. This is because increasing $P_m$ advances the machine 2 angle and therefore does not let machine 1 advance w.r.t machine 2 for the same fault clearing time.

Furthermore, as the trend is linear, a linear approximation holds well for large parameter variations as compared to the previous result. In both the results, similar to the case for the SMIB system, the CCT sensitivity computed (shown in red) is tangential to the CCT vs parameter curve. This proves the validity of the proposed approach.

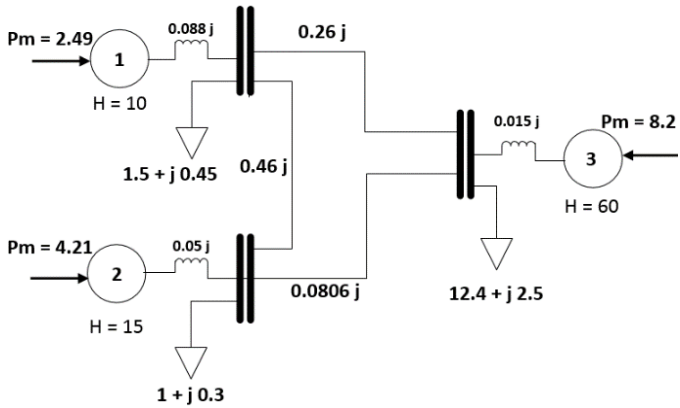

Figure 10 3-Machine System

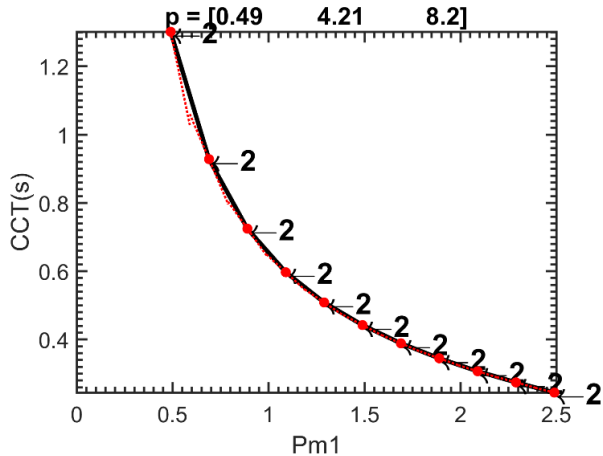

Figure 11 3 Machine System CCT vs $P_{m_1}$

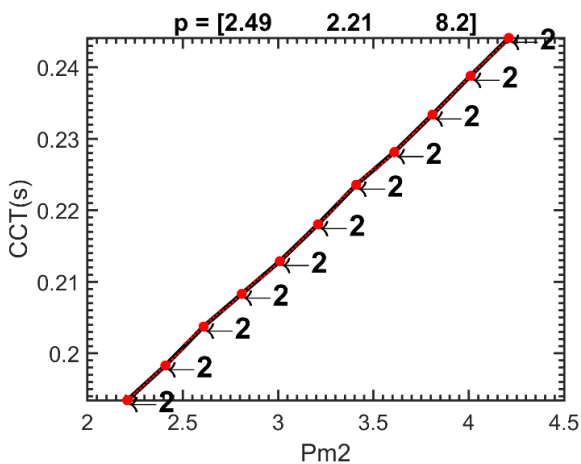

Figure 12 3 Machine System CCT vs $P_{m_2}$

## VIII. CONCLUSIONS AND FUTURE WORK

In this paper, given a critical fault-on and post-fault trajectory, we derived a formula for sensitivity of CCT to parameter variations for systems with inequality constraints. There are multiple mechanisms through which such systems become infeasible/unstable requiring a sensitivity formula derivation for each. A good application of this could be knowledge of the approximate impact of various system protection settings and operating conditions on changes in likelihood of undesirable tripping without using brute force methods.

It was observed that for constrained systems, the relevant portion of the constrained stability boundary might not be structurally stable under parameter variations unlike unconstrained systems, which are more robust. This situation occurred mainly due to involvement of some parameters being studied only limited to the system dynamics or feasibility constraints, but not both. We observed this when studying the variation of angle and frequency limits from the protection settings on CCT changes. This would require a more sophisticated approach to approximating CCT changes with parameter variations.

Lastly, we extended derivations obtained for ODE-governed systems to DAE-governed systems. The derivations for DAE systems would only be valid if the algebraic constraint results in a non-singular Jacobian. A case when this condition does not hold true is when (i) the load dynamics are not modeled, or (ii) along a trajectory passing through low voltage regions in state space and the maximum power transfer limit over some line(s) is reached, resulting in a loss of equilibrium of the algebraic system. These issues will be explored in our future work.